\documentclass[12pt]{amsart}

\usepackage{latexsym}
\usepackage[centertags]{amsmath}
\usepackage{amsfonts}
\usepackage{amssymb}
\usepackage{amsthm}
\usepackage{newlfont}
\usepackage[square]{natbib}

\def\mytopsep{3mm}

\newtheoremstyle{myplain}{\mytopsep}{\mytopsep}{\itshape}{0pt}{\bfseries}{.}{3mm}{}
\newtheoremstyle{mydefinition}{\mytopsep}{\mytopsep}{\normalfont}{0pt}{\bfseries}{.}{3mm}{}

\newtheoremstyle{myremark}{\mytopsep}{\mytopsep}{\normalfont}{0pt}{\bfseries}{.}{3mm}{}

\theoremstyle{myplain}

\newtheorem{thm}{Theorem}[section]

\newtheorem{lem}[thm]{Lemma}
\newtheorem{prop}[thm]{Proposition}

\theoremstyle{mydefinition}

\theoremstyle{myremark}
\newtheorem{rem}[thm]{Remark}
\newtheorem{exa}[thm]{Example}

\allowdisplaybreaks[1]

\makeatletter\@addtoreset{equation}{section}\makeatother

\def\sy{\mathfrak{S}}

\def\supp{\mathop{\mbox{supp}}}
\def\ord{\mathrm{ord}}

\def\NN{\mathbb{N}}

\def\xx{x^{-1}}
\def\yy{y^{-1}}
\def\CC{\mathbb{C}}
\def\ZZ{\mathbb{Z}}

\def\ct{\mathop{\mathrm{CT}}}

\def\pt{\mathop{\mathrm{PT}}}
\def\nt{\mathop{\mathrm{NT}}}
\def\CT{\mathop{\mathrm{CT}}}

\def\bous{Bousquet-M\'{e}lou}

\renewcommand{\ll}{\langle\!\langle}
\renewcommand{\gg}{\rangle\!\rangle}

\begin{document}

\title{Counting Lattice Paths By Gessel Pairs}

\author{Guoce Xin}
\address{Department
of Mathematics\\
Brandeis University\\
Waltham MA 02454-9110} \email{guoce.xin@gmail.com}

\date{September 14, 2004}
\begin{abstract}
We count a large class of lattice paths by using factorizations of
free monoids. Besides the classical lattice paths counting
problems related to Catalan numbers, we give a new approach to the
problem of counting walks on the slit plane (walks avoid a half
line) that was first solved by \bous\ and Schaeffer. We also solve
a problem about walks in the half plane avoiding a half line by
subsequently applying the factorizations of two different Gessel
pairs, giving a generalization of a result of \bous.
\end{abstract}
\maketitle


\section{Introduction}

 Ira Gessel (1980) \cite{ira} connected the factorization of formal Laurent series
 with that of lattice paths and solved a large class of lattice path
 enumeration problems, in which the unique factorization lemma (Lemma \ref{l-unifactor} below)
plays an important role.

After about two decades, \bous\ and Schaeffer \cite{bous}
remarkably solved the counting problem of walks on the slit plane:
lattice paths that start at $(0,0)$ with steps in a finite subset
$\sy$ of $\ZZ^2$ and never hit $(-k,0)$ for any nonnegative
integer $k$ after the starting point. Additionally they showed
that the complete generating functions of such walks are algebraic
in some models, and gave a surprising combinatorial interpretation
of Catalan numbers. In solving problems of walks on the slit plane
of \cite{bous,bouso}, the unique factorization lemma again plays
an important role.

This coincidence strongly suggests the existence of a
factorization of lattice paths that applies directly to walks on
the slit plane. We generalize Gessel's factorization of lattice
paths to that of a Gessel pair: a free monoid together with a
homomorphism to $\ZZ$. It turns out that such factorizations can
be repeatedly applied by choosing different Gessel pairs. This
technique yields a new approach to walks on the slit plane, and a
solution to the counting problems of lattice paths avoiding a half
plane and a half line.

We introduce the unique factorization lemma in Section 2, and the
concept of Gessel pairs in Section 3. In Section 4 we study some
explicit examples of Gessel pairs including walks on the slit
plane and walks one the half plane avoiding a half line. We shall
see that Theorem \ref{t-lagrange1} is a basic computational tool.

\section{The Unique Factorization Lemma}
The unique factorization lemma was developed in \cite{ira} for the
ring $K[[x,y/x]]$  and rediscovered in \cite{bous} for the ring
$K[x,y,\xx,\yy][[t]]$, where $K$ is a field. It will be rephrased
for the field of iterated Laurent series so that we can deal with
a larger class of problems.

Let $K$ be a field and let $x_1,\dots ,x_n$ be a set of variables.
The {\em field of iterated Laurent series} $K \ll x_1,\dots
,x_n\gg$ is inductively defined to be $K \ll x_1,\dots
,x_{n-1}\gg((x_n))$, with $K \ll x_1 \gg$ being the field of
Laurent series $K((x_1))$. An iterated Laurent series is first
regarded as a Laurent series in $x_n$, then a Laurent series in
$x_{n-1}$, and so on.

A fundamental structure theorem \cite{xinthesis} gives an overview
of iterated Laurent series:  {an iterated Laurent series} is a
formal series that has a well-ordered support, where $\ZZ^n$ is
ordered reverse lexicographically, and the \emph{support} of a
formal Laurent series is
$$\supp\ \sum_{(i_1,\dots ,i_n)\in \ZZ^n}
a_{i_1,\dots ,i_n} x_1^{i_1} \cdots x_n^{i_n} := \{\, (i_1,\dots ,
i_n)\mid a_{i_1,\dots ,i_n} \ne 0 \,\}.
$$
For an iterated Laurent series $f$, we can then define its order
$\ord(f)$ to be the minimum of its support, and its initial term
to be the term with the minimum order. We say that $f$ has a
positive order if $\ord(f)>\ord(1)=(0,\dots ,0)$.

Consequently, the following linear operators are well defined in
$\CC \ll x_1,\dots ,x_n\gg $:
$$\ct_x \sum_{n\in \ZZ} b_n x^n = b_0,\qquad \pt_x \sum_{n\in \ZZ} b_n x^n
= \sum_{n\ge 0} b_n x^{n},\qquad \nt_x \sum_{n\in \ZZ} b_n x^n =
\sum_{n<0} b_n x^{n},$$ where $x$ is one of the variables, and
$b_n$ may involve the other variables. We say that $f(x)$ is $\pt$
in $x$ if $f$ contains only nonnegative powers in $x$. Obviously
if $f(x)$ is $\pt$ in $x$ then $\ct_x f(x)=f(0)$. We will see that
these operators have combinatorial meanings in lattice path
enumeration.

The following two propositions follows from a general theory of
Malcev-Neumann series. See \cite[Theorem 3.1.7 and Proposition
3.2.6]{xinthesis}.

\begin{prop}[Composition Law]
Suppose that $F$ belongs to $K\ll x_1,\dots ,x_n\gg$ and
$\ord(F)>\ord(1)$. Then if $b_i\in K$ for all $i$,
$$\sum_{i=0}^\infty b_i f^i$$
is well defined and belong to $K\ll x_1,\dots ,x_n\gg$, in the
sense that all of its coefficients are finite sum of nonzero
elements in $K$.
\end{prop}

\begin{prop}[Generalized Composition Law]
Let $f$ be the initial term of $F\in K\ll x_1,\dots ,x_n\gg$. For
any $\Phi(x_1,\dots, x_n)$ belongs to $K\ll x_1,\dots ,x_n\gg$,
and any fixed $i$, $\Phi\mid _{x_i=F}$ is well defined if and only
if $\Phi\mid_{x_i=f}$ is well defined.
\end{prop}

By the composition law, if $\ord(F)>\ord(1)$, then $\log(1+F)$ is
well defined. The generalized composition law is useful in the
application of the kernel method. The other manipulation we will
use for iterated Laurent series are given as follows.

For any fixed variable $x$ and an iterated Laurent series $f(x)$,
we have that $f(x) $ can be uniquely written as $f_1(x)+f_2(x)$
with $f_1$ containing only nonnegative powers in $x$ and $f_2$
containing only negative powers in $x$. Clearly, $f_1(x)=\pt_x
f(x)$ and $f_2(x)=\nt_x f(x)$.

The unique factorization lemma follows from the above fact through
taking a logarithm.
\begin{lem}[Unique Factorization Lemma]\label{l-unifactor}
Let $h$ be an element of $K\ll x_1,\dots, x_n\gg$ with initial
term $1$. Then for each $i$, $h$ has a unique factorization in
$K\ll x_1,\dots ,x_n\gg$ such that $h=h_-h_0h_+$, where except for
their initial terms, which are $1$, $h_-$ contains only negative
powers in $x_i$, $h_0$ is independent of $x_i$, and $h_+$ contains
only positive powers in $x_i$.
\end{lem}

Theorem \ref{t-lagrange1} below is a generalization of the
Lagrange inversion formula. It plays an important role in the
proof of a conjecture about walks on the slit plane \cite{xin}.

\begin{thm}\label{t-lagrange1}
Let $G(x,t), F(x,t) \in K[[x,t]]$. If  $G(x,0)$ can be written as \\
$ax+\text{higher terms}$, with $a\ne 0$, then
\begin{equation}
\CT_{x} \frac{x}{G(x,t)} F(x,t)= \left.
\frac{F(x,t)}{\displaystyle{\partial\over
\partial x} G(x,t)} \right|_{x=X},
\end{equation}
where $X=X(t)$ is the unique element in $tK[[t]]$ such that
$G(X,t)=0$.
\end{thm}

\section{Factorization for Gessel Pairs}

By a {\em monoid} we mean a set $M$, equipped with a
multiplication which is associative and has a unit element $1$. An
element in $M$ is said to be a {\em prime} if it does not have a
nontrivial factorization. We say that $M$ is a {\em free monoid}
if every element in $M$ can be uniquely factored as a product of
primes in $M$.

We are going to present factorizations of free monoids with
respect to their homomorphisms to $\ZZ$. Our objects will be
mainly lattice paths: a path $\sigma$ in $\ZZ^2$ is a finite
sequence of lattice points $(a_0,b_0),\ldots , (a_n,b_n)$ in
$\ZZ^2$, in which we call $(a_0,b_0)$ the starting point,
$(a_n,b_n)$ the ending point, $(a_{i}-a_{i-1},b_i-b_{i-1})$ the
steps of $\sigma$, and $n$ the length of $\sigma$.

In this paper, the starting point of a path is always $(0,0)$. The
theory for other starting points is similar.

Given two paths $\sigma_1$ and $\sigma_2$, we define their product
$\sigma_1 \sigma_2$ to be the path whose steps are those of
$\sigma_1$ followed by those of $\sigma_2$. Thus the empty path
$\varepsilon$ is the unit. If $\pi=\sigma_1 \sigma_2$, then we
call $\sigma_1$ a {\em head} of $\pi$, and $\sigma_2$ a {\em tail}
of $\pi$.

Let $\sy$ be a finite subset of $\ZZ^2$. We are interested in
paths all of whose steps lie in $\sy$. Denote by $\sy^*$ the set
of all such paths. Then $\sy ^*$ is a free monoid and the primes
are elements of $\sy$. The weight of a step $(a,b)\in \sy $ is
defined to be $\Gamma((a,b))=x^ay^bt$, and the weight of a path
$\sigma =s_1\cdots s_n$ is defined to be $\Gamma(\sigma)=
\Gamma(s_1)\cdots \Gamma(s_n)$. It is easy to see that for any two
paths $\sigma_1$ and $\sigma_2$, we have
$\Gamma(\sigma_1\sigma_2)=\Gamma(\sigma_1)\Gamma( \sigma_2)$. If
$P$ is a subset of $\sy ^*$, then  we define
\begin{align}\label{e-lattice-p}
\Gamma(P)=\sum_{\sigma\in P}\Gamma(\sigma)=\sum_{n\ge 0}
\sum_{i,j\in \ZZ} a_{i,j}(n) x^iy^j t^n,
\end{align}
where $a_{i,j}(n)$ is the number of paths in $P$ of length $n$
that end at $(i,j)$. We also call $\Gamma(P)$ the generating
function of $P$ with respect to the ending points and the lengths.

In the special case that $P$ is the whole set $\sy ^*$, we have
$$\Gamma(\sy ^*)=\sum_{n\ge 0} (\Gamma(\sy ))^n=(1-\Gamma(\sy ))^{-1},$$
since each term in $(\Gamma(\sy ))^n$ corresponds to a path of $n$
steps.

The above equation is interpreted as an identity in
$\CC[x,\xx,y,\yy][[t]]$, which can be embedded into the field of
iterated Laurent series $\CC\ll x,y,t\gg$. In fact, we can relax
the condition on $\sy$ to a well-ordered subset of $\ZZ$, and the
composition law will guarantee the existence of $\Gamma(\sy^*)$
and hence $\Gamma(P)$.


Since $\sy $ is uniquely determined by $\Gamma(\sy )$, sometime we
give $\Gamma(\sy )$ instead of $\sy $. Some operators on $\CC\ll
x,y \gg[[t]]$ have simple combinatorial interpretations. Let $P$
be a subset of $\sy ^*$ with generating function given by
\eqref{e-lattice-p}.
\begin{enumerate}
\item The generating function for those paths in $P$
that end on the line $y=0$ is given by $\CT_y \Gamma(P)$.

\item The generating function for those paths in $P$
that end above the line $y=-1$ is given by $\pt_y \Gamma(P)$.

\item The generating function for those paths in $P$ that
end below the line $y=0$ is given by $\nt_y \Gamma(P)$.

\end{enumerate}
Similar properties hold for $x$.

\def\cxyt{\CC\ll x,y\gg[[t]]}

Now suppose that $H$ is a set of paths with steps in $\sy $ and
that $H$ is a free monoid. Then for any $\sigma \in H$ with its
factorization into primes $\sigma=h_1h_2\cdots h_m$, we say that
$h_1h_2\cdots h_i$ is an $H$-{\em head} of $\sigma$ for
$i=0,1,\ldots ,m$. If we let $\mathcal{P}$ be the set of primes in
$H$, then $\Gamma(H)=1/(1-\Gamma(\mathcal{P}))$.

For example, as we have described, $\sy ^*$ is a free monoid; the
set of all paths in $\sy ^*$ that end on the $x$-axis is a free
monoid, whose primes are those paths that return to the $x$-axis
only at the end point; the set of all paths in $\sy ^*$ that end
at $(k,0)$ for some $k\ge 0$ is a free monoid, whose primes are
those paths that only return the nonnegative half of the $x$-axis
at the end point.

Let $\rho$ be a map from $H$ to $\ZZ$. We say that $\rho$ is a
{\em homomorphism} from $H$ to $\ZZ$ if $\rho(\epsilon)=0$ and for
all $\sigma_1,\sigma_2\in H$, $\rho(\sigma_1 \sigma_2)
=\rho(\sigma_1)+\rho(\sigma_2)$. The $\rho$ value of a path
$\sigma$ is $\rho(\sigma)$.

If $H$ is a free monoid, then any map from $H$ to $\ZZ$ defined on
the primes of $H$ induces a homomorphism. If in addition, $H$ is a
subset of $\sy ^*$, then the natural map to the end point of a
path is a homomorphism from $H$ to $\ZZ^2$. Therefore, any
homomorphism from $\ZZ^2$ to $\ZZ$ induces a homomorphism from $H$
to $\ZZ$ through that natural map. The following two homomorphisms
are useful. Define $\rho_x(\sigma)$ to be  the $x$ coordinate of
the ending point of $\sigma$. Then $\rho_x$ is clearly a
homomorphism. Similarly we can define $\rho_y$.

If $H$ is a free monoid, and $\rho$ is a
 homomorphism from $H$ to $\ZZ$, then we call $(H,\rho)$
a \emph{Gessel pair}. For a Gessel pair $(H,\rho)$, we define:

A {\em minus-path} is either the empty path or a path whose $\rho$
value is negative and less than the $\rho$ values of all the other
$H$-heads.

A {\em zero-path} is a path with $\rho$ value $0$ and all of whose
$H$-heads have nonnegative $\rho$ values.

A {\em plus-path} is a path all of whose $H$-heads (except
$\epsilon$) have positive $\rho$ values.

For  a Gessel pair $(H,\rho)$, we denote by $H_ -$, $H_0$, and
$H_+$ respectively to be the sets of minus-, zero-, and plus-paths
in $H$. Note that the empty path, but no other path, belongs to
all three classes. The path $h_1h_2\cdots h_n$, where $h_i\in H$,
is a minus-path if and only if $h_nh_{n-1}\cdots h_1$ is a
plus-path; thus the theories of minus- and plus-paths are
identical.

\begin{lem}\label{l-pathfactor}
Let $(H,\rho)$ be a Gessel pair, and let $\pi$ be a path in $H$.
Then $\pi$ has a unique factorization $\pi_- \pi_0 \pi_+$, where
$\pi_-$ is a minus-path, $\pi_0$ is a zero-path, and $\pi_+$ is a
plus-path.
\end{lem}
\begin{proof}
Let $a$ be the smallest among all the $\rho$ values of the $H$
heads of $\pi$. Let $\pi_-$ be the shortest $H$-head of $\pi$
whose $\rho$ value equals $a$. Then if $\pi=\pi_- \sigma$, let
$\pi_-\pi_0$ be the longest $H$-head of $\pi$ whose $\rho$ value
equals $a$, and let $\pi_+$ be the rest of $\sigma$. It is easy to
see that this factorization satisfies the required conditions.

To see that it is unique, let $\tau_-\tau_0\tau_+$ be another
factorization of $\pi$. By definition, any $H$-head of $\tau_0
\tau_+$ has a nonnegative $\rho$ value. So the minimum $\rho$
value among all of the $H$-heads of $\pi$ is achieved in $\pi_-$.
By definition, it equals $\rho(\tau_-)$ and is unique in $\tau_-$.
Therefore, $\rho(\tau_-)=a$ and $\tau_-=\pi_-$ by the selection of
$\pi_-$. The reasons for $\pi_0=\tau_0$ and $\pi_+=\tau_+$ are
similar.

\end{proof}

\begin{prop}\label{p-pathfactor}
If $(H,\rho)$ is a Gessel pair, then $H_-$, $H_0$, and $H_+$ are
all free monoids. The map from $H$ to $H_-\times H_0\times H_+$
defined by $\pi\to(\pi_-,\pi_0,\pi_+)$ is a bijection.
\end{prop}
\begin{proof}
By Lemma \ref{l-pathfactor}, the map defined by
$\pi\to(\pi_-,\pi_0,\pi_+)$ is clearly a bijection. Now we show
that $H_-$, $H_0$, and $H_+$ are all free monoids.

It is easy to see that they are monoids. We only show that $H_-$
is free. The other parts are similar. Let $\mathcal{P}$ be the
subset of $H_-$ such that $\sigma \in \mathcal{P}$ if and only if
 $\rho(\sigma)$ is negative and every other
$H$-head of $\sigma$ has nonnegative $\rho$ value. We claim that
$\mathcal{P}$ is the set of primes in $H_-$.

Clearly any $\sigma\in \mathcal{P}$ cannot be factored as the
product of two nontrivial elements in $H_-$. Now let $\pi\in H_-$.
In order to factor $\pi$ into factors in $\mathcal{P}$, we find
the shortest $H$-head of $\pi$ that has negative $\rho$ value, and
denote it by $\sigma_1$. Then $\pi$ is factored as
$\pi=\sigma_1\pi'$ for some $\pi'$ in $H$. From the definition of
minus-path, $\rho(\sigma_1)$ is either less than $\rho(\pi)$, in
which case $\pi'$ is clearly in $ H_-$, or
$\rho(\sigma_1)=\rho(\pi)$, in which case $\pi'$ has to be the
unit and $\pi=\sigma_1$ is in $\mathcal{P}$. So we can inductively
obtain a factorization of $\pi$ into elements in $\mathcal{P}$.

The uniqueness of this factorization is clear.
\end{proof}

In a Gessel pair $(H,\rho)$, the weight of an element $\pi\in H$
is defined to be $\Gamma(\pi)z^{\rho(\pi)}$, where $z$ is a new
variable. When $H$ is also a subset of $\sy ^*$ and we are
considering the Gessel pair $(H,\rho_x)$, the power in  $z$ is
always the same as the power in $x$ for any $\pi$ in $H$. So we
can replace $z$ by $1$ and let $x$ play the same role as $z$.
Since the factorization in $H$ is with respect to $\rho$, the
factorization of generating function is with respect to $z$.

\begin{thm}\label{t-pathfactor}
For any Gessel pair $(H,\rho)$, we have
$\Gamma(H_-)=[\Gamma(H)]_-$, $\Gamma(H_0)=[\Gamma(H)]_0,$ and
$\Gamma(H_+)=[\Gamma(H)]_+$.
\end{thm}
\begin{proof}
From Proposition \ref{p-pathfactor}, it follows that
$\Gamma(H)=\Gamma(H_-)\Gamma(H_0)\Gamma(H_+)$. Clearly except $1$,
which is the weight of the empty path, $\Gamma(H_-)$ contains only
negative powers in $z$, $\Gamma(H_0)$ is independent of $z$, and
$\Gamma(H_+)$ contains only positive power in $z$. The theorem
then follows from the unique Factorization Lemma with respect to
$z$.
\end{proof}

Gessel \cite{ira} gives many interesting examples involving
lattice paths on the plane. We introduce the most classical
example as the following:
\begin{exa}
Let $\sy $ be $\{\, (1,r),(1,-1)\,\}$ with $r\ge 1$, and $H=\sy
^*$. Consider the Gessel pair $(H,\rho_y)$.
\end{exa}
Note that in this case the length of a path equals the $x$
coordinate of its end point. Replacing $x$ by $1$ will not lose
any information.

Clearly we have
$$\Gamma(H)=\Gamma(\sy ^*)=\frac{1}{1-t(y^r+1/y)}.$$

We see that $H_+$ is the set of paths in $\sy ^*$ that never go
below level $1$ after the starting point. The set $H_0$ contains
all paths in $\sy ^*$ that end on level $0$ and never go below
level $0$. When $r=1$, these are Dyck paths.

To compute $\Gamma(H_0):=F(t)$, we let $Y(t)$ be the unique
positive root of $y-t(1+y^{r+1})$. It is not hard to show that
$F(t)=Y(t)/t$ and hence $F(t)=1+t^{r+1}F(t)^{r+1}$. Therefore
$F(t)$ equals the generating function of complete $r+1$-ary trees.

\begin{exa}
Let $\sy $ be $\{\, (1,1), (1,-1) \,\}$, and let $H=\sy ^*$. Let
$\rho$ be determined by $\rho(1,1)=r$ and $\rho(1,-1)=-1$.
\end{exa}
It is easy to see that this example is isomorphic to the previous
one.

\begin{exa}\label{ex-4-hphl}
In general if $H=\sy ^*$, then $(H,\rho_y)$ is a Gessel pair.
\end{exa}
We see that $H_+$ is the set of paths in $\sy ^*$ that never go
below the line $y=1$ after the starting point.

If we let $J=H_+$, then $J$ is also a free monoid. The primes of
$J$ are paths that start at $(0,0)$, end at some positive level
$d$, and never hit level $d-1$ or lower.

The set $H_0$ contains all paths in $\sy ^*$ that end on the line
$y=0$, and never go below the line $y=0$. In other words, $H_0$
contains all paths in $\sy ^*$ that stays in the upper half plane
and end on the $x$-axis.

If we let $J=H_0$, then $(J,\rho_x)$ is a Gessel pair. The set
$J_+$ contains all paths in $J$ that avoid the half line
$\mathcal{H}=\{\, (-k,0); k\in \NN\,\}$ after the starting point.
This is the same as walks on the half plane avoiding the half line
in \citep{bouso}.

The set $J_0$ contains all paths in $J $ that end at $(0,0)$ and
never touch the half line $\mathcal{H}$ except $(0,0)$.

\vspace{3mm} {\em Walks on the slit plane} are paths that start at
$(0,0)$ with steps in $\sy$ and never hit the half line
$\mathcal{H}$ after the starting point. In solving counting
problem of walks on the slit plane \cite{bous,bouso}, it is
crucial to obtain the following functional equation
\eqref{e-slit-fac}, which will be explained combinatorially in
Example \ref{ex-4-slit}.

\begin{equation} \label{e-slit-fac}
S_0(x,t) \frac{1}{1-B(\xx,t)} =S_x(x;t)=\ct_y
\frac{1}{1-\Gamma(\sy)},
\end{equation}
where $B(\xx,t)$ is the generating function of paths that start at
$(0,0)$, and only hit $(-k,0)$ for some $k\ge 0$ at the end point;
$S_0(x,t)$ is the generating function of walks on the slit plane
that end on the line $y=0$; $S_x(x;t)$ is the generating function
of \emph{bilateral walks} \cite{bouso}: paths in $\sy^*$ that end
on the $x$-axis.

After obtaining equation \eqref{e-slit-fac}, we can check that
$S_x(x;0)=1$, $S_0(x,0)=1$, and $ B(\xx,0)=1$, and that except for
$1$, $S_0(x,t)$ contains only positive powers in $x$,
$(1-B(\xx,t))^{-1}$ contains only negative powers in $x$. Thus the
unique factorization lemma applies, and we obtain the following
remarkable result in \cite{bouso}, which says that the $S_0(x,t)$,
$B(\xx,t)$, and the complete generating function for walks on the
slit plane $S(x,y;t)$ can be theoretically computed. In practice,
computing them is not a easy task. Only special cases have been
thoroughly studied.

\begin{thm}[\bous]\label{t-4-slitg}
Let $\sy$ be a well-ordered subset in $\ZZ^2$. Using notation as
above, we have:
\begin{align}
S_0(x,t) &=\left(S_x(x,t)\right)_+, \\
\frac{1}{1-B(\xx,t)} &= \left(S_x(x,t)\right)_0\left(S_x(x,t)\right)_-,\\
S(x,y;t)&=
\frac{1}{(1-\Gamma(\sy))\left(S_x(x,t)\right)_0\left(S_x(x,t)\right)_-}.\label{e-4-sxyt}
\end{align}
\end{thm}

Walks on the slit plane can be counted by a factorization of
Gessel pair.
\begin{exa}\label{ex-4-slit}
For any $\sy $, let $H$ be the set of paths that end on the
$x$-axis. Then $(H,\rho_x)$ is a Gessel pair.
\end{exa}
The set $H_+$ contains all paths that end on the $x$ axis and
never hit the half line $\mathcal{H}=\{\, (-k,0)\mid k\ge 0\,\}$
after the starting point. This is exactly the {walks on the slit
plane} that end on the $x$-axis.

The set $H_0$, which was called the set of loops in \citep{bouso},
consists of all paths that end at $(0,0)$, and never touch
$(-k,0)$ for $k=1,2,\ldots$.

The combinatorial explanation of equation \eqref{e-slit-fac} is as
follows. The set $H_-H_0$ is a free monoid. It contains all paths
that end at $(-k,0)$ for some $k\ge 0$. Its primes are all paths
that hit $(-k,0)$ only once at the end point. Clearly, the primes
are are counted by $B(\xx;t)$. So we have
$$\Gamma(H_-H_0)=\frac{1}{1-B(\xx;t)}, \quad \text{ and } \Gamma(H_+)=S_0(x,t).$$
Equation \eqref{e-slit-fac} then follows.

\begin{rem}
Note that equation \eqref{e-slit-fac} can also be explained
combinatorially by using the cycle lemma as in \citep{bouso}.
\end{rem}

\begin{exa}\label{ex-4-hphl1}
For any $\sy $, let $H$ be the set of paths that end on the
$x$-axis and never go below the line $y=-d$ for some given $d>0$.
Then it is easy to check that $(H,\rho_x)$ is a Gessel pair.
\end{exa}

The set $H_+$ contains all paths that end on the $x$-axis, and
never hit the half line $\mathcal{H}$ after the starting point,
and never go below the line $y=-d$.

The set $H_0$ can be similarly described.

\begin{exa}\label{ex-4hphl2}
For any $\sy $, let $H$ be the set of paths that end on the
$x$-axis and never go below the line $y=-d$ and never go above the
line $y=f+1$ for some given positive integers $d$ and $e$. Then it
is easy to see that $(H,\rho_x)$ is a Gessel pair.
\end{exa}
This example is similar to the previous one.

\section{Explicit Examples\label{s-4-exa}}

We will discuss two explicit examples that were proposed in
\citep{bouso}. Taking walks on the slit plane as an example, we
see that $\log S_0(x,t)=\pt_x \log S_x(x;t)$. Now if $\log
S_0(x,t)$ has the form $b(t)x^p+\text{higher degree terms}$, then
so does $S_0(x,t)-1$.

\begin{prop}[Proposition 4, \citep{bouso}]\label{p-sp0}
Let $p$ be the smallest positive integer such that there is a walk
on the slit plane with respect to a finite set $\sy$ that ends at
$(p,0)$. Then the generating function for such walks ending at
$(p,0)$ is $D$-finite and is given by
$$S_{p,0}(t)=[x^p] \log S_x(x;t).$$
\end{prop}
\bous\ shows in addition that $S_{k,0}$ is $D$-finite for every
$k$.

Our task is to find a formula for $\log \Gamma(H_+)$ for a given
algebraic $\Gamma(H)$ as described in last section. The idea is as
follows. Let $P(x,y,t)$ be a polynomial and let $Y(x;t)$ be the
unique root of positive order of $y-tP(x,y,t)$ as a polynomial in
$y$. The problem will be reduced to finding the unique
factorization of a rational function $Q(x,Y(t),t)$ with respect to
$x$. We are especially interested in $[x^p]Q_+(x,Y(t),t)$ for a
certain integer $p$, which is $D$-finite by Proposition
\ref{p-sp0}. This generating function can be obtained if we can
get a nice form for $\frac{\partial}{\partial t}\log Q(x,Y(t),t)$.
Our approach to finding such a nice form is to do all the
computation implicitly. It is best illustrated by examples.

\begin{exa}
Let $\sy$ be the set $\{ \, (1,0),(-1,0),(0,2)(0,-1)\, \}$, or
equivalently, $\Gamma(\sy)=t(x+\xx+y^2+\yy)$. \bous\ proposed in
\citep{bouso} the problem of solving walks on the slit plane in
this model, or even replacing the $2$ by a general positive
integer $q$.
\end{exa}
Our method works for general $q$, but so far we have found a
reasonable formula only for $q=2$. We have:

\begin{prop}
The number of walks on the slit plane, of length $N$, ending at
$(1,0)$, and with steps in  $\{ \, (1,0),(-1,0),(0,2),(0,-1)\, \}$
equals
\begin{multline}\label{e-4-exay2s10}
a_{1,0}(N)= \binom{N}{\frac{N-1}2}+\sum_{n=1}^{\lfloor N/3\rfloor
}\frac{3^{3n-1}}{n2^{2n}}\binom{N-1}{3n-1}\binom{N-3n}{\frac{N-3n}{2}}+\\
\sum_{n,m,k}\frac{3^{3m+2}}{nN2^{2m+2}}
\binom{n}{k,2k+1,n-3k-1}\binom{N-n}{3m+2}
\binom{N-3m-3k-3}{\frac{N-3m-3k-4}{2}},
\end{multline}
where $\binom{A}{B+1/2}$ is interpreted as $0$ for all integers
$A,B$, and the second sum ranges over all $n,m,k$ such that $1\le
n\le N$, $0\le m\le \frac{N-n-2}{3} $, and $0\le k\le
\frac{n-1}{3}$.
\end{prop}

\begin{proof}
We proceed by computing $S_x(x;t)$. Let $b=x+\xx$. Then
$\Gamma(\sy)=t(b+y^2+\yy)$. Applying Theorem \ref{t-lagrange1}, we
get
\begin{align*}
S_x(x;t)&=\ct_y \frac{y}{y-t(y^3+by+1)}=\frac{1}{1-tb-3tY^2},
\end{align*}
where $Y=Y(t)=Y(b,t)=Y(x,t)$ is the unique root of positive order
of the denominator for $y$. More precisely, $Y$ is the unique
power series in $t$ with constant term $0$ that satisfies
\begin{align}
\label{e-4-exay2} Y(t)-t(Y(t)^3+bY(t)+1)=0.
\end{align}
Using the Lagrange inversion formula we get
\begin{align}\label{4-exay2Y}
Y(t)=\sum_{n\ge 1} \sum_{k=0}^{\left\lfloor
\frac{n-1}{3}\right\rfloor} \binom{n}{k,2k+1,n-3k-1} b^{n-3k-1}
t^n.
\end{align}

We can compute $\log S_x(x;t)$ explicitly in order to obtain $\log
S_0(x,t)$. We have
\begin{align}
\frac{\partial }{\partial t}\log S_x(x;t)&=\frac{b+3 \left( Y
\left( t \right) \right) ^{2}+6tY \left( t \right) {
\frac{\partial }{\partial t}}Y \left( t
\right) }{1-tb-3tY(t)^2}\nonumber \\
&=\frac{b-t{b}^{2}+3\, \left( Y \left( t \right)  \right) ^{2}-3\,
\left( Y
 \left( t \right)  \right) ^{4}t+6\,tY \left( t \right)
}{(1-tb-3tY(t)^2)^2}\label{e-4-exay2-lns}
\end{align}
where
$$ \frac{\partial }{\partial t}Y \left( t \right) = \frac{1+bY(t)+Y(t)^3}{1-tb-3tY(t)^2}$$ is
determined implicitly by equation \eqref{e-4-exay2}.

Since $Y(t)$ satisfies \eqref{e-4-exay2}, we can rewrite
\eqref{e-4-exay2-lns} as $C_0+C_1Y(t)+C_2Y(t)^2$, where $C_i$ are
rational functions of $b$ and $t$. These can be found by Maple,
and we get
\begin{align}\label{e-4-exay2-lnsf}
\frac{\partial }{\partial t}\log S_x(x;t)= \frac{\left(
4\,{b}^{3}+27 \right) {t}^{2}-8\,t{b}^{2}+4\,b }{4(1-bt)^3-27t^3}+
\frac{9tY(t)}{4(1-bt)^3-27t^3}.
\end{align}
Now we need to integrate to get $\log S_x(x;t)$. The first term
has a simple form:
$$ \int \frac{\left( 4\,{b}^{3}+27 \right)
{t}^{2}-8\,t{b}^{2}+4\,b }{4(1-bt)^3-27t^3}dt = \log
\left(4(1-bt)^3-27t^3\right)^{-1/3}+C,
$$
where $C$ is independent of $t$. After some manipulation, we get
\begin{align*} \int \frac{\left( 4\,{b}^{3}+27
\right) {t}^{2}-8\,t{b}^{2}+4\,b }{4(1-bt)^3-27t^3}dt
=\log\frac{1}{1-bt}+\sum_{N\ge 1} \sum_{n=1}^{\lfloor
N/3\rfloor}\frac{3^{3n-1}}{n2^{2n}}\binom{N-1}{3n-1}b^{N-3n}t^N+C.
\end{align*}

For the second term, we have
\begin{align*}
\frac{9t}{4(1-bt)^3-27t^3} &= \frac{9t}{4(1-bt)^3}
\frac{1}{1-27t^3/(4(1-bt)^3)}.
\end{align*}
After some manipulation, we get
\begin{align*}
\frac{9t}{4(1-bt)^3-27t^3} &= \sum_{m\ge 0}
\frac{3^{3m+2}}{2^{2m+2}} \sum_{r\ge 0 } \binom{3m+r+2}{3m+2} b^r
t^{3m+r+1}.
\end{align*}
Thus together with the expansion of $Y(t)$ given by
\eqref{4-exay2Y}, we obtain
\begin{multline*}
    \int \frac{9tY(t)}{4(1-bt)^3-27t^3}dt =C+\sum_{N\ge
1}\sum_{n=1}^N\sum_{m=0}^{\left\lfloor \frac{N-n-2}{3}
\right\rfloor }
\sum_{k=0}^{\left\lfloor \frac{n-1}{3} \right\rfloor}\\
    \frac{3^{3m+2}}{nN2^{2m+2}}
\binom{n}{k,2k+1,n-3k-1}\binom{N-n}{3m+2}b^{N-3m-3k-3} t^N.
\end{multline*}
Since $S_x(x,0)=1$, it is easy to check that the sum of the two
constants $C$ must be $0$.

Note that the powers in $b$ is always nonnegative. It is easy to
separate the negative powers and positive powers in $x$ of
$b^M=(x+\xx)^M$ for every nonnegative integer $M$. Thus we can
obtain a formula for $\log S_0(x,t)$. In particular, from the
formulas $[x] (x+\xx)^M=\binom{M}{\frac{M-1}{2}}$ and
$S_{1,0}(t)=[x]\log S_x(x;t)$, we get \eqref{e-4-exay2s10}.
\end{proof}

\begin{exa}
We consider walks on the half plane avoiding a half line; more
precisely, walks that never touch the half line $\mathcal{H}$ and
never hit a point $(i,j) $ with $j<0$. This is a continuation of
Example \ref{ex-4-hphl}. We denote by $HS(x,y;t)$ the generating
function for such paths.
\end{exa}

It turns out that this case is simpler than the previous one. We
obtain the following result, which includes \citep[Proposition
25]{bouso} as a special case.
\begin{thm}
For any well-ordered set $\sy$, let $p$ be the smallest positive
number such that there is an  $\sy$-path end at $(p,0)$. Then the
number of walks on the half plane avoiding the half line that end
at $(p,0)$ and are of length $n$ is equal to $1/n$  times the
number of $\sy$-paths that end at $(p,0)$ and are of length $n$.
\end{thm}

\begin{proof}
We use the notation of Example \ref{ex-4-hphl}. From the Gessel
pair $(\sy^*, \rho_y)$, we have
$\Gamma(H_0)=\left(\Gamma(\sy^*)\right)_0$ and
\begin{align*}
\log \Gamma(H_0) =\ct_y \log \Gamma(\sy^*)=\ct_y \log
\frac{1}{1-\Gamma(\sy)}.
\end{align*}

Now let $J=H_0$ and consider the Gessel pair $(J, \rho_x)$. Then
$$ \log \Gamma(J_0J_+)= \pt_x  \log \Gamma(J). $$
In particular, we have
\begin{align*}
[x^p] \Gamma(J_+)&= [x^p] \log \Gamma(J)=[x^p] \log
\Gamma(H_0)=[x^p] \ct_y \log \Gamma(\sy^*).
\end{align*}
Therefore,
$$[x^p t^n] \Gamma(J_+) = [x^py^0t^n] \frac{1}{n}\Gamma(\sy)^n. $$
This prove the theorem.
\end{proof}

\vspace{3mm} {\bf Acknowledgment.} The author is very grateful to
Ira Gessel and Mireille \bous.

\bibliographystyle{amsplain}

\providecommand{\bysame}{\leavevmode\hbox
to3em{\hrulefill}\thinspace}
\providecommand{\MR}{\relax\ifhmode\unskip\space\fi MR }
\providecommand{\MRhref}[2]{%
  \href{http://www.ams.org/mathscinet-getitem?mr=#1}{#2}
} \providecommand{\href}[2]{#2}

\end{document}